\documentclass[dvips, 12pt, a4paper]{article}
\usepackage{latexsym,amsmath,amssymb}
\numberwithin{equation}{section}

\begin{document}

\begin{center}
{\Large \textbf{Statistical approximation properties of
$(p,q)$-Szasz-Mirakjan Kantorovich operators}}

\bigskip

\textbf{Bhausaheb R. Sontakke} and \textbf{Amjad Shaikh}
\bigskip

Department of Mathematics, Pratishthan Mahavidyalaya, Paithan - 431 107 Dist. Aurangabad, (M.S.) India\\[0pt]
Department of Mathematics, Poona College of Arts, Science and Commerce, Camp, Pune- 411001, (M.S.), India\\[0pt]

brsontakke@rediffmail.com ; amjatshaikh@gmail.com \\[0pt]

\bigskip

\textbf{Abstract}
\end{center}

{\footnotesize {The main aim of this study is to introduce statistical
approximation properties of $(p,q)$-Szasz Mirakjan Kantorovich operators with the help of the
Korovkin type statistical approximation theorem. Rates of statistical
convergence by means of the modulus of continuity and the Lipschitz type
maximal function are also established.
}}\newline



\bigskip

\section{ Introduction}

\parindent=8mm \parindent=8mm   Recently, Mursaleen et al \cite{mka1} applied $(p,q)$-calculus in
approximation theory and introduced the first $(p,q)$-analogue of Bernstein
operators based on $(p,q)$-integers. Motivated by the work of Mursaleen et al \cite{mka1}, the idea of $(p,q)$-calculus and its importance.
Very recently, Khalid et al. has given a very nice application in computer-aided geometric design and applied these Bernstein basis  for construction of $(p,q)$-B$\acute{e}$zier curves and surfaces based on  $(p,q)$-integers which is further generalization of $q$-B$\acute{e}$zier curves and surfaces \cite{wcq,khalid,khalid1}. For similar works based on $(p,q)$-integers, one can refer \cite{acar1,acar2,acar3,zmn,mur8,mka3,mnak1,ma1,mah,wafi,cai,ali}.\\

\parindent=8mm It was  S.N. Bernstein \cite{brn} in 1912, who first introduced his famous operators $%
B_{n}: $ $C[0,1]\rightarrow C[0,1]$ defined for any $n\in \mathbb{N}$
and for any function $f\in C[0,1]$
\begin{equation}\label{e1}
B_{n}(f;x)=\sum\limits_{k=0}^{n}\left(
\begin{array}{c}
n \\
k%
\end{array}%
\right) x^{k}(1-x)^{n-k}f\biggl{(}\frac{k}{n}\biggl{)},~~x\in \lbrack 0,1].
\end{equation}
 and named it Bernstein polynomials to prove the Weierstrass
theorem \cite{pp}.
Later it was found that Bernstein polynomials possess many remarkable properties and has various applications in areas such as approximation theory \cite{pp}, numerical analysis,
computer-aided geometric design, and solutions of differential equations due to its fine properties of approximation \cite{hp}.

In  computer aided geometric design (CAGD), Bernstein polynomials and its variants are used in order to preserve the shape of the curves or surfaces. One of the most important curve in CAGD \cite{thomas} is the classical B$\acute{e}$zier curve \cite{Bezier} constructed with the help of Bernstein basis functions.

 The  Szasz-Mirakjan operators\cite{GM, OSz}  have an important role in the approximation theory, and their approximation
properties have been investigated by many researchers. The Kantorovich type of the
Szasz Mirakjan operators was defined as
\begin{equation}
K_{n}(f;x)= ne^{-nx}\sum_{k=0}^{\infty}\frac{(nx)^k}{k!}\int_{k/n}^{(k+1)/n }%
f(t)dt,\,\,\,\,\,\,\,\,\, f\in C_{\gamma}[0,\infty )
\end{equation}

\begin{equation*}
where   \ C_{\gamma}[0,\infty ) = \{f\in C_{\gamma}[0,\infty ): \mid f(t)\mid \leq M {(1 + t)}^\gamma
\end{equation*}
In\cite{arg}, Aral and Gupta defined q-type generalization of Szasz Mirakjan operators as follows:
\begin{equation}
S_{n}^{q}(f;q,x)= E_{q}(-[n]_{q}\frac{x}{b_{n}})\sum_{k=0}^{\infty}f(\frac{[k]_{q}b_{n}}{[n]_{q}})\frac{([n]_{q}x)^{k}}{[k]_{q}!(b_{n})^{k}}
\end{equation}

where $ f\in C[0,\infty ), q \in (0, 1), 0\leq x<\frac{b_{n}}{1-q_{n}}, b_{n}$ is a sequence of positive numbers such that $\lim_{n\rightarrow \infty }b_{n}= \infty \,\,\, and \,\,\, E_{q}(x)= \sum_{n=0}^{\infty}q^{\frac{n(n-1)}{2}}\frac{x^{n}}{[n]_{q}!}$\newline

Let $C_{B}[0,\infty )$ be the space of all bounded and continuous functions
on $[0,\infty )$ .Then $C_{B}[0,\infty )$ is a normed linear space with $%
\Vert f\Vert _{C_{B}}=\sup\limits_{x\geq 0}\mid f(x)\mid $. Let $w$ be a
function of the type of modulous of continuity. The principal properties of
the function are the following:\newline
$(i)~~$ $w$ is a nonnegative increasing function on $[0,\infty )$,\newline
$(ii)~~$ $\lim\limits_{\delta \rightarrow 0}~w(\delta )=0$.\newline

Let $H_{w}$ be the space of all-real valued functions $f$ defined
on $[0,\infty )$ satisfying the following condition:\newline
\begin{equation*}
\mid f(x)-f(y)\mid \leq w(\mid x-y \mid )\newline
\end{equation*}%
for any $x,y\geq 0.$\newline

In \cite{Gdj}, Gadjiev and Caker proved the Korovkin type theorem which gives the conditions of the convergence of the sequence of linear positive operators to find function in $H_{w}$.\newline
Currently, more useful connections of Korovkin type approximation theory, not only with classical approximation theory, but also other branches of mathematics were given by Altomare and Campiti in\cite{campiti}.\newline
Now we recall the following theorem which was given by gadjiev and Cakar:\newline

\textbf{Theorem 1.1}(\cite{Gdj}). Let $(A_{n})$ be the sequence of linear positive operators from $H_{w}$ into $C_{B}[0,\infty )$ satisfying three conditions
\begin{equation*}
\lim_{n\rightarrow \infty }\Vert A_{n}(t^{\nu}; x)-x^{\nu}\Vert _{c_{B}}=0, \,\,\,\, \nu = 0, 1, 2.
\end{equation*}
Then for any function ${f}\in$ $H_{w}$
\begin{equation*}
\lim_{n\rightarrow \infty }\Vert A_{n}(f)- f\Vert _{c_{B}}=0
\end{equation*}

Let us give rudiments of $$(p, q)$$-calculus.\newline
 For each nonnegative integer $n$, the $(p,q)$%
-integer $[n]_{p,q}$ is defined by\newline
\begin{equation*}
 [n]_{p,q}:= \frac{(p^{n}-q^{n})}{(p-q)}, \,\,\,   n = 0,1,2,...,\,\,\,\,  0<q<p\leq1
\end{equation*}
whereas q-integers are given by
\begin{equation*}
 [n]_{q}:= \frac{(1-q^{n})}{(1-q)}, \,\,\,   n = 0,1,2,...,\,\,\,\,  0<q\leq1
\end{equation*}
and the $(p,q)$-binomial coefficients are defined by
\begin{equation*}
\left[
\begin{array}{cc}
n &\\
k &
\end{array}%
\right] _{p,q}:=\frac{[n]_{p,q}!}{[k]_{p,q}![n-k]_{p,q}!}
\end{equation*}%

By some simple calculation, we have the following relation
\begin{equation*}
q^{k}[n-k+1]_{p,q}=[n+1]_{p,q}-p^{n-k+1}[k]_{p,q}
\end{equation*}
For details on $q$-calculus and $(p,q)$-calculus, one is referred to\cite {vp,sad}\newline
Recently in \cite{makh}, Kantorovich Variant of $(p,q)$-Szász-Mirakjan Operators has been studied and defined as follows\newline
$For f\in C[0,\infty ), 0<q<p\leq1 $ and each positive integer n
\begin{equation}
K_{n}(f,p,q; x)= [n]_{p,q} \sum_{k=0}^{\infty}p^{-k}q^{k}s_{n,k}(p,q;x)\int_{[k]_{p,q}/q^{k-1}[n]_{p,q}}^{[k+1]_{p,q}/q^{k}[n]_{p,q} }%
f(t)d_{p,q}t
\end{equation}
where f is a nondecreasing function.\newline
In\cite{makh} M. Mursaleen et al. obtained the uniform approximation of these operators to the function $f\in H_w$ as follows \newline

\textbf{Theorem 1.2}(\cite{makh}). Let $q = q_{n} , p = p_{n}$  be a sequence satisfying $0<q_{n}<p_{n}\leq1 $ and let $ q_{n}\rightarrow1\,\,, p_{n}\rightarrow1 \,\,\,as\,\,\, n\rightarrow\infty$. If $K_{n}$ is defined by(1.3), then for any $f\in H_w$ \newline
\begin{equation*}
\lim_{n\rightarrow \infty }\Vert K_{n}(f)- f\Vert _{c_{B}}=0
\end{equation*}%

On the other hand the concept of statistical convergence was introduced by Fast \cite{fast} in the year
1950 and in recent times it has become an active area of research. The concept
of the limit of a sequence has been generalized to a statistical limit through the
natural density of a set K of positive integers, defined as
\begin{equation*}
\delta(K) = \lim_{n\rightarrow \infty} \frac{1}{n}\,\,(k\leq n \,\,for\,\, k \in K)
\end{equation*}
provided this limit exists\cite{zman}. We say that the sequence $x = (x_{n}) $ statistically
converges to a number L. if for each $\epsilon > 0$, the density of the set $\{k\in
\mathbb{N}:|x_{k}-L|\geq \epsilon \}$.We denote it by  $st-\lim\limits_{k}x_{k}=L$. It is easily seen that every convergent
sequence is statistically convergent but not inversely.\newline
Statistical convergence was used in approximation theory by Gadjiev and Orhan \cite{gadorh}. They proved the statistically Korovkin type theorem for the linesr positive operators as follows

\textbf{Theorem 1.3} (\cite{gadorh}). If the sequence of positive linear operators $A_{n}: C[a, b]\rightarrow C[a, b]$
satisfies the conditions
\begin{equation*}
st-\lim_{n\rightarrow \infty }\Vert A_{n}(e_{\nu}; .)-e_{\nu} \Vert _{c_{[a,b]}}=0
\end{equation*}%
$with \,\,\,e_{\nu} = t^{\nu} \,\,\,for \,\,\,\nu = 0, 1, 2 $  \newline
then  for any function\,\, $f\in C[a, b],$ we have
\begin{equation*}
st-\lim_{n\rightarrow \infty }\Vert A_{n}(f; .)- f \Vert _{c_{[a,b]}}=0
\end{equation*}

The main aim of this paper is to study the operators defined by M.Mursaleen et.al \cite{makh} and obtain statistically approximation properties of the operator with the help of the Korovkin type theorem proved by A.D. gadjiev and C. Orhan \cite{gadorh} and estimate the rate of statistically convergence of the sequence of the operator to the function f.\newline
we had obtained the ordinary result for the uniform convergence of the bivariate $(p, q)$--Szasz-Mirakjan operators in \cite{makh}. In this study, the similar results will be used while proving the theorems.\newline
\section{Main Results:}

Let us give the following theorem:\newline
\textbf{Theorem 2.1 :} Let $(A_{n})$ be the sequence of linear positive operators from $H_{w} $ into $ C_{B}(R_{+}) $ satisfying three conditions
\begin{equation*}
st-\lim_{n\rightarrow \infty }\Vert A_{n}(t^{\nu}: x)- x^{\nu} \Vert _{c_{B}}=0,\,\,\,\,\,\,\,\,\,\,\,\,\,\,\,\,\,\,\,\,\, \nu = 0, 1, 2
\end{equation*}
then for any $ f\in H_{w} $
\begin{equation*}
st-\lim_{n\rightarrow \infty }\Vert A_{n}(f; .)- f \Vert _{c_{B}}=0
\end{equation*}
Notice that this theorem is one variable case of the Duman and Erkus theorem \cite{erkusd}\newline
To obtain the statistically convergence of the operators above, we need the following three lemmas as given in \cite{makh}\newline

\textbf {Lemma 2.2}~~~For $ x \geq 0$, \,\,\, $ 0< q <p \leq1 $\newline

$(i) \,\,\,  K_{n}(1, p, q; x) = 1$

$(ii) \,\,\,  K_{n}(t, p, q; x) = \frac{1}{q}x + \frac{1}{[2]_{p, q}[n]_{p,q}}$

$(iii) \,\,\,  K_{n}(t^{2}, p, q; x) = \frac{p}{q^{3}}x^{2} + (\frac{p + [2]_{p,q}}{q[3]_{p,q}[n]_{p,q}} + \frac{1}{q^{2}[n]_{p,q}})x + \frac{1}{[3]_{p,q}[n]_{p,q}^{2}}$\newline

\section{Korovkin type statistical approximation properties}

The main aim of this paper is to obtain the Korovkin type statistical approximation properties of operators defined in(1.3), with the help of Theorem (1.1).\newline
Now, we consider a sequence $ p = p_{n}$, $q = q_{n}$  satisfying the following expression
\begin{equation}
 st-\lim\limits_{n}q_{n}= 1,~~~ st-\lim\limits_{n}p_{n}= 1,~~~ 0<q_{n}<p_{n}\leq1
\end{equation}
\textbf{Theorem 3.1} \,\,\,\, Let $(K_{n})$ be the sequence of the operators (1.3) and the sequence $q = q_{n}$,\,\,\, $p = p_{n}$ satisfies (1.3) for $0<q_{n}<p_{n}\leq1$ then for any $f\in H_{w}$
\begin{equation*}
st-\lim_{n\rightarrow \infty }\Vert K_{n}(f; q_{n};p_{n}; .)- f \Vert _{c_{B}}=0
\end{equation*}
Proof:  In the light of Theorem 2.1, it is sufficient to prove the followings:
\begin{equation*}
st-\lim_{n\rightarrow \infty }\Vert K_{n}(t^{\nu}: x)- x^{\nu} \Vert _{c_{B}}=0,\,\,\,\, \nu = 0, 1, 2
\end{equation*}
From Lemma , the first condition of above equation is easily obtained for $\nu = 0.$\newline
for $\nu = 1$
\begin{equation*}
\Vert K_{n}(t;q_{n};p_{n}; x )- x \Vert _{c_{B}}= \Vert \frac{1}{q_{n}}x + \frac{1}{[2]_{p_{n}, q_{n}}[n]_{p_{n},q_{n}}}- x \Vert
\leq \,\,\mid (\frac{1}{q_{n}} - 1) + \frac{1}{[2]_{p_{n}, q_{n}}[n]_{p_{n},q_{n}}} \mid
\end{equation*}
Now for a given $\epsilon > 0,$ we define following sets\newline
\begin{equation*}
U = \{n : \Vert K_{n}(t;q_{n};p_{n}; x )- x \Vert \geq \epsilon \}
\end{equation*}
\begin{equation*}
U_{1} = \{n : 1 -\frac{1}{q_{n}} \geq \epsilon \}, \,\,\,\,\,\, U_{2} = \{n : \frac{1}{[2]_{p_{n}, q_{n}}[n]_{p_{n},q_{n}}} \geq \epsilon \}
\end{equation*}
It is obvious that $ U \subset U_{1}\cup  U_{2} $
Therefore it can be written as\newline
\begin{equation*}
\delta\{k \leq n :\Vert K_{n}(t; q_{n};p_{n}; x )- x \Vert \geq \epsilon \} \leq \delta\{k \leq n:1 -\frac{1}{q_{n}}  \geq \epsilon \} + \delta\{k \leq n :\frac{1}{[2]_{p_{n}, q_{n}}[n]_{p_{n},q_{n}}} \geq \epsilon \}
\end{equation*}
By using (3.1), it is clear that \newline
\begin{equation*}
st-\lim_{n\rightarrow \infty }(1 -\frac{1}{q_{n}}) = 0\,\,\,and \,\,\,st-\lim_{n\rightarrow \infty }(\frac{1}{[2]_{p_{n}, q_{n}}[n]_{p_{n},q_{n}}}) = 0
\end{equation*}
\begin{equation*}
so \,\,\, \delta\{k \leq n:1 -\frac{1}{q_{n}}  \geq \epsilon \} = 0 \,\,\,and \,\,\, \delta\{k \leq n :\frac{1}{[2]_{p_{n}, q_{n}}[n]_{p_{n},q_{n}}} \geq \epsilon \} = 0
\end{equation*}
\begin{equation*}
st-\lim_{n\rightarrow \infty }\Vert K_{n}(t; q_{n};p_{n}; .)- x \Vert _{c_{B}}=0
\end{equation*}
Lastly for $\nu = 2$ we have
\begin{equation*}
\Vert K_{n}(t^{2}; q_{n};p_{n}; x)- x^{2} \Vert
= sup_{x \geq 0}\bigg\{ x^{2}(\frac{p_{n}}{q_{n}^{3}} - 1) + x(\frac{p_{n} + [2]_{p_{n},q_{n}}}{q_{n}[3]_{p_{n},q_{n}}[n]_{p_{n},q_{n}}} + \frac{1}{q_{n}^{2}[n]_{p_{n},q_{n}}}) + \frac{1}{[3]_{p_{n},q_{n}}[n]_{p_{n},q_{n}}^{2}}\bigg\}
\end{equation*}

~~~~~~~~~~~~~~~~~~~~~~~$\leq|\frac{p_{n}}{q_{n}^{3}} - 1| +|\frac{p_{n} + [2]_{p_{n},q_{n}}}{q_{n}[3]_{p_{n},q_{n}}[n]_{p_{n},q_{n}}} + \frac{1}{q_{n}^{2}[n]_{p_{n},q_{n}}}| + |\frac{1}{[3]_{p_{n},q_{n}}[n]_{p_{n},q_{n}}^{2}}|$
\newline

~~~~~~~~~~~~~~~~~~~~~~~$=\frac{p_{n}}{q_{n}^{3}} - 1 + \frac{p_{n} + [2]_{p_{n},q_{n}}}{q_{n}[3]_{p_{n},q_{n}}[n]_{p_{n},q_{n}}} + \frac{1}{q_{n}^{2}[n]_{p_{n},q_{n}}} + \frac{1}{[3]_{p_{n},q_{n}}[n]_{p_{n},q_{n}}^{2}}$\newline

If we choose $ \alpha_{n} = \frac{p_{n}}{q_{n}^{3}} - 1 $, $ \beta_{n} = \frac{p_{n} + [2]_{p_{n},q_{n}}}{q_{n}[3]_{p_{n},q_{n}}[n]_{p_{n},q_{n}}} + \frac{1}{q_{n}^{2}[n]_{p_{n},q_{n}}} $,  $\gamma_{n} = \frac{1}{[3]_{p_{n},q_{n}}[n]_{p_{n},q_{n}}^{2}} $ \newline

 then one can write
\begin{equation}
  st-\lim_{n\rightarrow \infty }\alpha_{n} = st - \lim_{n\rightarrow \infty }\beta_{n} = st - \lim_{n\rightarrow \infty }\gamma_{n} = 0
\end{equation}
 by (3.1)
 Now given $\epsilon > 0$, we define the following four sets
\begin{equation*}
U = \{n : \Vert K_{n}(t^{2};q_{n};p_{n}; x )- x^{2} \Vert \geq \epsilon \},
\end{equation*}
\begin{equation*}
U_{1} = \{n:\alpha_{n}\geq \frac{\varepsilon}{3}\},\,\,\,U_{2} = \{n:\beta_{n}\geq \frac{\varepsilon}{3}\}\,\,\,U_{3} = \{n:\gamma_{n}\geq \frac{\varepsilon}{3}\}
\end{equation*}
It is obvious that $ U\subseteq U_{1}\cup U_{2} \cup U_{3}.$  Then we obtain
\newline

$\delta\{k \leq n :\Vert K_{n}(t^{2}; q_{n};p_{n}; x )- x^{2} \Vert \geq \epsilon \}\leq \delta\{k \leq n:\alpha_{n}\geq \frac{\epsilon}{3}\} + \delta\{k \leq n :\beta_{n}\geq \frac{\epsilon}{3}\}$
\newline

~~~~~~~~~~~~~~~~~~~~~~~~~~~~~~~~~~~~~~~~~~~~~~~~~$~+~\delta\{k \leq n :\gamma_{n}\geq \frac{\epsilon}{3}\}$

So the right hand side of the inequalities is zero by (3.2), then
\begin{equation*}
st-\lim_{n\rightarrow \infty }\Vert K_{n}(t^{2}; q_{n};p_{n}; .)- x^{2}\Vert _{c_{B}}=0
\end{equation*}

holds. Hence the proof follows from theorem (2.1).\\

\section{ Rates of statistical convergence }

\parindent=8mmIn this section, we give the rates of statistical convergence
of the operator (1.3) by means of modulus of continuity and Lipschitz type
maximal functions. The modulus of continuity for the functions $f\in H_{w} $ is defined as\newline
\begin{equation*}
w(f;\delta )=\sup\limits_{x,t\geq 0,\,\,\mid t-x\mid < \delta }\mid f(t)-f(x)\mid
\end{equation*}%
\newline
where $w(f;\delta )~\text{for}~\delta >0 $
satisfies the following conditions: for every $ f\in H_{w}$
\newline
$(i)$ $\lim\limits_{\delta \rightarrow 0 }w(f;\delta)=0$\newline
$(ii)$ $\mid f(t)-f(x){\mid }\leq w(f;\delta )\left(\frac{\mid
t-x\mid}{\delta}+1\right)~~~~~~~~~~~~~~~~~~~~~~~~~~~~~~~~~~~~~~~~~~~~~~~~~~~~(4.1)$ \newline

\textbf{Theorem 4.1.} Let the sequence $q=(q_{n})$, $p=(p_{n})$ satisfies the condition in $%
(3.1)$ and $0<q_{n}< p_{n}\leq1$. Then we have\newline
\begin{equation*}
\mid K_{n}(f;q_{n};p_{n}; x)-f(x)\mid \leq 2 w(f;\sqrt{\delta _{n}(x)})
\end{equation*}%
where
\newline
$ \delta _{n}(x) =  x^{2}(\frac{p_{n}}{q_{n}^{3}} - \frac{2}{q_{n}} + 1)
 + x(\frac{p_{n} +  [2]_{p_{n},q_{n}}}{q_{n}[3]_{p_{n},q_{n}}[n]_{p_{n},q_{n}}} + \frac{1}{q_{n}^{2}[n]_{p_{n},q_{n}}} - \frac{2}{[2]_{p_{n},q_{n}}[n]_{p_{n},q_{n}}})+\frac{1}{[3]_{p_{n},q_{n}}[n]_{p_{n},q_{n}}^{2}}$~~~~(4.2)\newline

\textbf{Proof.}\newline
Since $\mid K_{n}(f;q_{n};p_{n}; x)-f(x)\mid \leq K_{n}({\mid f(t)-f(x)\mid };q_{n}; p_{n}; x),$ by $(4.1)$ we get\newline
\begin{equation*}
\mid K_{n}(f;q_{n};p_{n}; x)-f(x)\mid \leq w(f;\delta )~~\left(\{K_{n}(1;q_{n};p_{n}; x)+\frac{1}{\delta_{n}}K_{n}({\mid t-x\mid };q_{n};p_{n};x)\}\right).
\end{equation*}%
\newline
Using Cauchy-Schwartz inequality, we have
\begin{align*}
\mid K_{n}(f;q_{n};p_{n}; x)-f(x)\mid &\leq w(f;\delta _{n})\left(1+%
\frac{1}{\delta _{n}}[(K_{n}({\ t-x)^{2}}%
;q_{n};p_{n};x)]^{\frac{1}{2}}{[K_{n}(1;q_{n};p_{n};x)]}^{\frac{1}{2}}\right)\\
&\leq w(f;\delta _{n})\bigg\{1 +\frac{x^{2}}{\delta _{n}}\bigg(\frac{p_{n}}{q_{n}^{3}} - \frac{2}{q_{n}} + 1\bigg)\\
&+\frac{x}{\delta _{n}}\bigg(\frac{p_{n} +  [2]_{p_{n},q_{n}}}{q_{n}[3]_{p_{n},q_{n}}[n]_{p_{n},q_{n}}}
 + \frac{1}{q_{n}^{2}[n]_{p_{n},q_{n}}} - \frac{2}{[2]_{p_{n},q_{n}}[n]_{p_{n},q_{n}}}\bigg)\\
 &+\frac{1}{\delta _{n}}\bigg(\frac{1}{[3]_{p_{n},q_{n}}[n]_{p_{n},q_{n}}^{2}}\bigg)\bigg\}
\end{align*}

By choosing $\delta _{n}$ as in $(4.2)$, we get the desired result.

This completes the proof of the theorem.\newline

Note that in condition $(3.1)$,%
\begin{equation*}
st-\lim_{n\rightarrow \infty }\delta _{n}=0.
\end{equation*}%
By $(4.1)$ we have
\begin{equation*}
st-\lim_{n\rightarrow \infty }w(f;\delta _{n})=0,
\end{equation*}%
which gives us the pointwise rate of statistical convergence of the operator
$K_{n}(f;q_{n};p_{n}x)$ to $f(x).$\newline
Now we will give an estimate concerning the rate of approximation by means
of Lipschitz type maximal functions. In \cite{lenze},B. Lenze introduced a Lipschitz
type maximal function as
\begin{equation*}
\tilde{f}_{\alpha }(x)=\sup\limits_{t>0,t\neq x}\frac{\mid f(t)-f(x)\mid }{{%
\mid {x-t}\mid }^{\alpha }}.
\end{equation*}%
\newline
In \cite{ardog}, the Lipschitz type maximal function space on $E\subset \lbrack
0,\infty )$ is defined as follows
\begin{equation*}
\tilde{W}_{\alpha }=\{f=\sup (1+x)^{\alpha }~\tilde{f}_{\alpha }(x)\leq M%
\frac{1}{(1+y)^{\alpha }};x\geq 0~{and}~y~\in E\},
\end{equation*}
where $f$ is bounded and continuous function on $[0,\infty )$, $M$ is a
positive constant and $0<\alpha \leq 1$. \newline
we denote by d(x, E), the distance between x and E, that is $ d(x,E)=inf \{\mid x - y \mid; y \in E\}$

\textbf{Theorem 4.2.} If $K_{n}$ be defined by $(1.3)$, then for all $%
f\in \tilde{W}_{\alpha ,E}$
\begin{equation*}
\mid K_{n}(f;q_{n};p_{n};x)-f(x)\mid \leq M\left(\delta _{n}^{\frac{\alpha }{2}}+ 2~(d(x,E))^{\alpha}\right)~~~~~~~~~~~~~~~~~~~~~~~~~~~~~~~~~~~~~~~~~~~~~~(4.3)
\end{equation*}%

where $ \delta_{n}(x) $ is defined in Theorem (4.1)

\textbf{Proof.} Let $x\geq 0,~(x,x_{0})\in \lbrack 0,\infty )\times E$. Then
we have\newline
\begin{equation*}
\mid f-f(x)\mid \leq \mid f-f(x_{0})\mid +\mid f(x_{0})-f(x)\mid .
\end{equation*}%
Since $K_{n}$ is a positive and linear operator, $f\in \tilde{W}%
_{\alpha ,E}$ and using the above inequality
\begin{equation*}
\mid K_{n}(f;q_{n};p_{n};x)-f(x)\mid \leq \mid K_{n}(\mid
f-f(x_{0})\mid ;q_{n};p_{n};x)-f(x)\mid +\mid f(x_{0})-f(x)\mid
K_{n}(1;q_{n};p_{n};x)
\end{equation*}%

~~~~~~~~~~~~~~~~~~~~~~~~~$\leq M \left( K_{n}({\mid t-x_{0}\mid }%
^{\alpha };q_{n};p_{n};x)+{\mid x-x_0 \mid}^{\alpha }K_{n}(1;q_{n};p_{n};x)\right)$
\newline
Therefore we have
\begin{equation*}
K_{n}\left( {\mid t-x_{0}\mid }^{\alpha
};q_{n};p_{n};x\right) \leq K_{n}({\mid t-x\mid }%
^{\alpha };q_{n};p_{n};x)+{\mid x-x_{0}\mid}^{\alpha }K_{n}(1;q_{n};p_{n};x).
\end{equation*}%
By using the Holder inequality with $p=\frac{2}{\alpha }$ and $q=\frac{2}{%
2-\alpha }$, we have
\newline
$K_{n}\left( \mid t-x\mid^{\alpha };q_{n};p_{n};x\right)\leq K_{n}{\left( {(t-x)}^{2};q_{n};p_{n};x\right)}^\frac{\alpha}{2}{(K_{n}(1;q_{n};p_{n};x))}^{\frac{2-\alpha }{2}}
+{\mid x-x_{0}\mid} ^{\alpha }K_{n}(1;q_{n};p_{n};x)$\newline

~~~~~~~~~~~~~~~~~~~~~~~~$=\delta _{n}^{\frac{\alpha }{2}}+{\mid x-x_{0}\mid} ^{\alpha}.$

This completes the proof of the theorem.\newline

\textbf{Remark 4.1 }If we take $E=[0,\infty )$ in Theorem (4.2), since $d(x,E)=0,$ then we obtain the following result:\newline
For every $f\in \tilde{W}_{\alpha ,[0,\infty )}$
\begin{equation*}
{\mid K_{n}(f;q_{n};p_{n}x)-f(x)\mid }\leq M\delta _{n}^{\frac{%
\alpha }{2}}.
\end{equation*}%
where $\delta_{n}$ is defined as in $(12).$\newline
\textbf{Remark 4.2} By using $(3.1)$, It is easy to verify that
\begin{equation*}
st-\lim_{n\rightarrow \infty }\delta _{n}=0.
\end{equation*}%
That is, the rate of statistical convergence of $(1.3)$ are estimated by
means of Lipschitz type maximal functions.

 \section{Construction of the bivariate operators}

In what follows we construct the bivariate extension of the operators (1.3).
We will introduce the statistical convergence of the operators to a function f and
investigate the statistical rate of convergence of these operators.\newline

 Let $ R_{+}^{2} = [0, \infty) \times [0, \infty) ,\,\,\, f: R_{+}^{2}\rightarrow R $ and $ 0< p_{n_{1}},q_{n_{1}};p_{n_{2}},q_{n_{2}}\leq1$ Then
we define the bivariate companion of the operators (1.3) as follows:\newline

$ K_{n_{1},n_{2}}(f,p_{n_{1}},q_{n_{1}};p_{n_{2}},q_{n_{2}}; x,y) = $
\begin{equation*}
 = [n_{1}]_{p_{n_{1}},q_{n_{1}}}[n_{2}]_{p_{n_{2}},q_{n_{2}}}
 \sum_{k_{1}=0}^{\infty}\sum_{k_{2}=0}^{\infty}p_{n_{1}}^{-k_{n_{1}}}q_{n_{1}}^{k_{n_{1}}}
 s_{n_{1},k_{1}}(p_{n_{1}}, q_{n_{1}}; x,y)p_{n_{2}}^{-k_{n_{2}}}q_{n_{2}}^{k_{n_{2}}}
 s_{n_{2},k_{2}}(p_{n_{2}}, q_{n_{2}}; x,y)
\end{equation*}
\begin{equation}
\times\int_{[k_{1}]_{p_{n_{1}},q_{n_{1}}}/q_{n_{1}}^{k_{1}-1}[n_{1}]_{p_{n_{1}},q_{n_{1}}}}^{[k_{1}+1]_{p_{n_{1}},q_{n_{1}}}/q_{n_{1}}^{k_{1}}[n_{1}]_{p_{n_{1}},q_{n_{1}}} }\int_{[k_{2}]_{p_{n_{2}},q_{n_{2}}}/q_{n_{2}}^{k_{2}-1}[n_{2}]_{p_{n_{2}},q_{n_{2}}}}^{[k_{2}+1]_{p_{n_{2}},q_{n_{2}}}/q_{n_{2}}^{k_{2}}[n_{2}]_{p_{n_{2}},q_{n_{2}}} }f(t)d_{p_{n_{1}},q_{n_{1}}}t~d_{p_{n_{2}},q_{n_{2}}}t
\end{equation}
\newline
For$ K = [0,\infty) × [0,\infty)$, the modulus of continuity for the bivariate case is defined as
\begin{equation*}
 w_{2}(f; \delta_{1}, \delta_{2})= sup{\mid f(u, v)- f(x, y)\mid : (u, v), (x, y) \in K and \mid u - x \mid \leq \delta_{1}, \mid v - y \mid \leq \delta_{2}}
\end{equation*}

 and  $ w_{2}(f; \delta_{1}, \delta_{2})$ satisfy the following condition
\begin{equation*}
 \mid f(u, v)- f(x, y)\mid \leq  w_{2}\left( f ;\mid u - x \mid, \mid v - y \mid \right)
\end{equation*}

for each $ f \in H_{w_{2}}$. Detailed study of modulus of continuity for the bivariate analogue one is
referred to \cite{anas}.\newline
The first Korovkin type theorem for the statistical approximation for the bivariate
analogue of linear positive operators defined in the space $H_{w_{2}}$ was obtained
by Erkus and Duman \cite{erkusd} which is as follows

\textbf{Theorem 5.1}\cite{erkusd}. Let $K_{n}$ be a sequence of positive linear operators from  $H_{w_{2}}$ into $C_{B}(K)$. Then for each $f\in H_{w_{2}}$,
\begin{equation*}
  st - \lim_{n\rightarrow \infty }\Vert K_{n}(f)- f\Vert =0
\end{equation*}
 is satisfied if the following holds:\newline
 \begin{equation*}
 st - \lim_{n\rightarrow \infty }\Vert K_{n}(f_{i})- f\Vert =0, \,\,\,\,\,\,\,\,\,\,\,\, i = 0,1,2,3
 \end{equation*}
 where

 \,\,\,\,$ f_{0}(u, v) = 1$, \,\,\,\,    $ f_{1}(u, v) = u $, \,\,\,    $ f_{2}(u, v) = v $, \,\,\,\,\,\,\,  $ f_{3}(u, v) = u^{2} + v^{2} $~~~~~~~~~~~~~~~~~~~~~~~~~~~~~~~~~~~(5.2)
 \newline

To study the statistical convergence of the bivariate operators, the following
\newline

\textbf{Lemma 5.2.} The bivariate operators defined above satisfy the followings: \newline
$(i)\,\,\, K_{n_{1},n_{2}}(f_{0},p_{n_{1}},q_{n_{1}};p_{n_{2}},q_{n_{2}}; x,y) = 1 $\newline
$(ii)\,\,\, K_{n_{1},n_{2}}(f_{1},p_{n_{1}},q_{n_{1}};p_{n_{2}},q_{n_{2}}; x,y) = \frac{1}{q_{n_{1}}}x + \frac{1}{[2]_{p_{n_{1}}, q_{n_{1}}}[n_{1}]_{p_{n_{1}},q_{n_{1}}}} $\newline
$(iii)\,\,\, K_{n_{1},n_{2}}(f_{1},p_{n_{1}},q_{n_{1}};p_{n_{2}},q_{n_{2}}; x,y) = \frac{1}{q_{n_{2}}}y + \frac{1}{[2]_{p_{n_{2}}, q_{n_{2}}}[n_{2}]_{p_{n_{2}},q_{n_{2}}}} $\newline
$(iv)\,\,\, K_{n_{1},n_{2}}(f_{1},p_{n_{1}},q_{n_{1}};p_{n_{2}},q_{n_{2}}; x,y) = x^{2}(\frac{p_{n_{1}}}{q_{n_{1}}^{3}} - 1) + x(\frac{p_{n_{1}} + [2]_{p_{n_{1}},q_{n_{1}}}}{q_{n_{1}}[3]_{p_{n_{1}},q_{n_{1}}}[n_{1}]_{p_{n_{1}},q_{n_{1}}}} + \frac{1}{q_{n_{1}}^{2}[n_{1}]_{p_{n_{1}},q_{n_{1}}}}) + \frac{1}{[3]_{p_{n_{1}},q_{n_{1}}}[n_{1}]_{p_{n_{1}},q_{n_{1}}}^{2}} + y^{2}(\frac{p_{n_{2}}}{q_{n_{2}}^{3}} - 1) + y(\frac{p_{n_{2}} + [2]_{p_{n_{2}},q_{n_{2}}}}{q_{n_{2}}[3]_{p_{n_{2}},q_{n_{2}}}[n_{2}]_{p_{n_{2}},q_{n_{2}}}} + \frac{1}{q_{n_{2}}^{2}[n_{2}]_{p_{n_{2}},q_{n_{2}}}}) + \frac{1}{[3]_{p_{n_{2}},q_{n_{2}}}[n_{2}]_{p_{n_{2}},q_{n_{2}}}^{2}} $\newline

Proof. Exploiting the proofs for the bivariate operators in \cite{ersan}, the above can be
easily established. So we skip the proof.
\newline Now, we consider a sequence $ p = p_{n_{1}}$, \,\,\,$ p = p_{n_{2}}$,\,\,\,$ q = q_{n_{1}}$, \,\,\,$ q = q_{n_{2}}$,  be statistically convergent to unity but not convergent in usual sense, so we can
write them for\newline
~~~~~~~~~~~~~~~~~~~~~~~~~$ 0 < p_{n_{1}}, p_{n_{2}}, q_{n_{1}}, q_{n_{2}} \leq 1$
\begin{equation*}
 st-\lim\limits_{n}q_{n_{1}}= st-\lim\limits_{n}p_{n_{1}}= st-\lim\limits_{n}q_{n_{2}}= st-\lim\limits_{n}p_{n_{2}}=1 ~~~~~~~~~~~~~~~~~~~~~~~~~~~~~~~~~~~~~~~~~~~~~(5.3)
\end{equation*}
Now under the condition in (5.3), let us show that the statistical convergence of the bivariate operator (5.1) with the help of the proof of Theorem (3.1).\newline

\textbf{Theorem 5.3.} Let $ p = (p_{n_{1}}) $, $p = (p_{n_{2}})$, $q = (q_{n_{1}})$ and $q = (q_{n_{2}})$ be the sequences
satisfy the conditions (4.3) and let $K_{n_{1}}, K_{n_{2}} $ be the sequence of linear positive oper-
ators from $H_{w_{2}}(R^{2}_{+})$ into $C_{B}(R_{+})$. Then for each $f\in H_{w_{2}}$,

\begin{equation*}
  st - \lim_{n_{1},n_{2}\rightarrow \infty }\Vert K_{n_{1},n_{2}}(f)- f\Vert = 0
\end{equation*}
Proof. With the aid of the Lemma (5.2), a proof similar to the proof of the
Theorem (3.1) can be easily obtained. So we shall omit the proof.\newline

 \section{Rates of convergence of the bivariate operators}

For any $f \in H_{w_{2}}(R^{2}_{+})$, the modulus of continuity of the bivariate analogue is defined as:
\begin{equation*}
 \tilde{w}(f; \delta_{1}, \delta_{2})= \sup\limits_{t, x \geq 0}\left\{\mid f(t, s)- f(x, y)\mid ;   \mid t-x\mid \leq\delta_{1}, \mid s-y\mid \leq\delta_{2}, \,\,\,(t, s), (x, y) \in R^{2}_{+}\right\}
\end{equation*}
For details of this sort of modulus, one is referred to \cite{anas}. Here $\tilde{w}(f; \delta_{1}, \delta_{2})$ satisfies the following conditions\newline

$(i)~~\tilde{w}(f; \delta_{1}, \delta_{2}) \rightarrow 0 \,\,if\,\, \delta_{1}\rightarrow 0, \,\,\,\delta_{2}\rightarrow 0 $

$(ii)~~\mid f(t, s)- f(x, y)\mid \leq \tilde{w}(f; \delta_{1}, \delta_{2})\left(\frac{\mid
t-x\mid}{\delta_{1}}+1\right)\left(\frac{\mid s-y\mid}{\delta_{2}}+1\right)~~~~~~~~~~~~~~~~~~~~~~~~~~~~~~~~~~~~~~~~~~~~~(6.1)$ \newline

Now in the following theorem we study the rate of statistical convergence of
the bivariate operators through modulus of continuity in $H_{w_{2}}$. \newline

\textbf{Theorem 6.1.} Let $p = (p_{n_{1}})$, $p = (p_{n_{2}})$, $q = (q_{n_{1}})$, $q = (q_{n_{2}})$ be four sequences
obeying conditions of (4.3). Then we have
\begin{equation*}
 \mid K_{n_{1},n_{2}}(f,p_{n_{1}},q_{n_{1}};p_{n_{2}},q_{n_{2}}; x,y)- f(x,y)\mid \leq 4 w\left(f;\sqrt{\delta _{n_{1}}(x)};\sqrt{\delta _{n_{2}}(y)}\right)
\end{equation*}

where\newline
$ \delta _{n_{1}}(x) =  x^{2}(\frac{p_{n_{1}}}{q_{n_{1}}^{3}} - \frac{2}{q_{n_{1}}} + 1)
 + x(\frac{p_{n_{1}} +  [2]_{p_{n_{1}},q_{n_{1}}}}{q_{n_{1}}[3]_{p_{n_{1}},q_{n_{1}}}[n_{1}]_{p_{n_{1}},q_{n_{1}}}} + \frac{1}{q_{n_{1}}^{2}[n_{1}]_{p_{n_{1}},q_{n_{1}}}} - \frac{2}{[2]_{p_{n_{1}},q_{n_{1}}}[n_{1}]_{p_{n_{1}},q_{n_{1}}}})$\newline

~~~~$+\frac{1}{[3]_{p_{n_{1}},q_{n_{1}}}[n_{1}]_{p_{n_{1}},q_{n_{1}}}^{2}}$\newline
$\delta _{n_{2}}(y) =  y^{2}(\frac{p_{n_{2}}}{q_{n_{2}}^{3}} - \frac{2}{q_{n_{2}}} + 1)
 + y(\frac{p_{n_{2}} +  [2]_{p_{n_{2}},q_{n_{2}}}}{q_{n_{2}}[3]_{p_{n_{2}},q_{n_{2}}}[n_{2}]_{p_{n_{2}},q_{n_{2}}}} + \frac{1}{q_{n_{2}}^{2}[n_{2}]_{p_{n_{2}},q_{n_{2}}}} - \frac{2}{[2]_{p_{n_{2}},q_{n_{2}}}[n_{2}]_{p_{n_{2}},q_{n_{2}}}})$\newline

~~~~$+\frac{1}{[3]_{p_{n_{2}},q_{n_{2}}}[n_{2}]_{p_{n_{2}},q_{n_{2}}}^{2}}$\newline

Proof: Using the property of the modulus above, we have\newline
\begin{equation*}
\mid K_{n_{1},n_{2}}(f,p_{n_{1}},q_{n_{1}};p_{n_{2}},q_{n_{2}}; x,y)- f(x,y)\mid \leq w(f; \delta_{n_{1}}, \delta_{n_{2}})\,\bigg(\{K_{n_{1},n_{2}}(f_{0};p_{n_{1}},q_{n_{1}};p_{n_{2}},q_{n_{2}}; x,y)
\end{equation*}
\begin{equation*}
+\frac{1}{\delta_{n_{1}}}K_{n_{1},n_{2}}({\mid t-x\mid };p_{n_{1}},q_{n_{1}};p_{n_{2}},q_{n_{2}}; x,y)\}\{K_{n_{1},n_{2}}(f_{0};p_{n_{1}},q_{n_{1}};p_{n_{2}},q_{n_{2}}; x,y)
\end{equation*}
\begin{equation*}
+\frac{1}{\delta_{n_{2}}}K_{n_{1},n_{2}}({\mid s-y\mid };p_{n_{1}},q_{n_{1}};p_{n_{2}},q_{n_{2}}; x,y)\}\bigg)
\end{equation*}
Applying the Cauchy-Schwartz inequality, we get
\begin{equation*}
K_{n_{1},n_{2}}({\mid t-x\mid };p_{n_{1}},q_{n_{1}};p_{n_{2}},q_{n_{2}}; x,y)\leq[(K_{n_{1},n_{2}}({\ t-x)^{2}}%
;p_{n_{1}},q_{n_{1}};p_{n_{2}},q_{n_{2}}; x,y)]^{\frac{1}{2}}
\end{equation*}
\,\,\,\,\,\,\,\,\,\,\,\,\,\,\,\,\,\,\,\,\,\,\,\,\,\,\,\,\,\,\~~~~~~~~~~~~~~~~~~~~~~~~~~~~~~~~~~~~~
$\times {[K_{n_{1},n_{2}}(f_{0};p_{n_{1}},q_{n_{1}};p_{n_{2}},q_{n_{2}}; x,y)]}^{\frac{1}{2}}$\newline
On substituting this in the above inequality, we get the proof of the theorem.

Now we shall study the statistical convergence of the bivariate operators
using Lipschitz type maximal functions.\newline
The Lipschitz type maximal function space on $ E \times E \subset R_{+} \times R_{+}$ is defined as follows

\begin{equation*}
\tilde{W}_{\alpha_{1},\alpha_{2}}E^{2}=\{f: \sup (1+t)^{\alpha_{1}}(1+s)^{\alpha_{2}}~\tilde{f}_{\alpha_{1},\alpha_{2}}(x,y)\leq M%
\frac{1}{(1+x)^{\alpha_{1}}}\frac{1}{(1+y)^{\alpha_{2}}};x,y\geq 0, (t,s)\in E^{2}\}
\end{equation*}
~~~~~~~~~~~~~~~~~~~~~~~~~~~~~~~~~~~~~~~~~~~~~~~~~~~~~~~~~~~~~~~~~~~~~~~~~~~~~~~~~~~~~~~~~~~~~~~~~~~~~~~~~~~~~~~~~~~~~~~~~~~~~~~~(6.2)
Where f is a bounded and continuous function on $R_{+},$ M is a positive constant
and $ 0\leq \alpha_{1}, \alpha_{2} \leq 1$ and $\tilde{f}_{\alpha_{1},\alpha_{2}}(x, y)$ is defined as follows:

\begin{equation*}
\tilde{f}_{\alpha_{1},\alpha_{2}}(x, y) = \sup\limits_{t, s\geq 0}\frac{\mid f(t,s)-f(x,y)\mid }{{%
\mid {t-x}\mid }^{\alpha_{1} }{\mid {s-y}\mid }^{\alpha_{2} }}.
\end{equation*}

\textbf{Theorem 6.2.} Let $ p = (p_{n_{1}}),~ p = (p_{n_{2}}), ~ q = (q_{n_{1}}), ~ q = (q_{n_{2}}) $ be four sequences
satisfying the conditions of (5.2). Then we have
\begin{equation*}
\mid K_{n_{1},n_{2}}(f;p_{n_{1}},q_{n_{1}};p_{n_{2}},q_{n_{2}}; x,y)-f(x,y)\mid \leq M_{p_{n_{1}},q_{n_{1}},p_{n_{2}},q_{n_{2}}}\bigg({\delta_{n_{1}}(x)^{\frac{\alpha_{1}}{2}}{\delta_{n_{2}}(y)^{\frac{\alpha_{2}}{2}}}(p_{n_{1}},q_{n_{1}},p_{n_{2}},q_{n_{2}})}
\end{equation*}
~~~~~~~~~~~~~~~~~~~~~~~~~~~~~~~~~~~~~~~~~~~~~~~~~~$ +~{\delta_{n_{1}}(x)^{\frac{\alpha_{1}}{2}}}d(y,E)^{\alpha_{2}} + {\delta_{n_{2}}(y)^{\frac{\alpha_{2}}{2}}}d(x,E)^{\alpha_{1}} +~2~d(x,E)^{\alpha_{1}}d(y,E)^{\alpha_{2}}\bigg) $
where \newline
$0 \leq \alpha_{1}, \alpha_{2} \leq 1 ~~and~~\delta_{n_{1}}, \delta_{n_{2}} $ are defined as in Theorem (6.1) and \newline
$ d(x,E) = inf \{|x-y| : y \in E\}.$\newline
Proof. For $ x, y \geq 0 $ and $(x_{1}, y_{1}) \in E \times E, $ we can write
\begin{equation*}
\mid f(t,s)-f(x,y)\mid \leq \mid f(t,s)-f(x_{0},y_{0})\mid +\mid f(x_{0},y_{0})-f(x,y)\mid .
\end{equation*}%
Applying the operator $K_{n_{1},n_{2}}$ to both sides of the above inequality and making
use of eqn (6.2), we have
\begin{equation*}
\mid K_{n_{1},n_{2}}(f;p_{n_{1}},q_{n_{1}};p_{n_{2}},q_{n_{2}}; x,y)-f(x,y)\mid \leq \mid K_{n_{1},n_{2}}(\mid
f(t,s)-f(x_{0},y_{0})\mid ;p_{n_{1}},q_{n_{1}};p_{n_{2}},q_{n_{2}}; x,y)
\end{equation*}
~~~~~~~~~~~~~~~~~~~~~~~~~~~~~~~~~~~~~~~~~~~~~~~~~~~~~~~~~~~~~~~$+\mid f(x_{0},y_{0})-f(x,y)\mid K_{n_{1},n_{2}}(f_{0};p_{n_{1}},q_{n_{1}};p_{n_{2}},q_{n_{2}}; x,y)$\newline

~~~~~~~~~~~~~~~~~~~~~~~~~~~~~~~~~~~~$\leq M K_{n_{1},n_{2}}({\mid t-x_{0}\mid }^{\alpha_{1}}{\mid s-y_{0}\mid }^{\alpha_{2}};p_{n_{1}},q_{n_{1}};p_{n_{2}},q_{n_{2}}; x,y)$\newline

~~~~~~~~~~~~~~~~~~~~~~~~~~~~~~~~~~~~$+M{\mid x-x_{0} \mid}^{\alpha_{1}}{\mid y-y_{0} \mid}^{\alpha_{2}}K_{n_{1},n_{2}}(f_{0};p_{n_{1}},q_{n_{1}};p_{n_{2}},q_{n_{2}}; x,y)$

~~~~~~~~~~~~~~~~~~~~~~~~~~~~~~~~~~~~~~~~~~~~~~~~~~~~~~~~~~~~~~~~~~~~~~~~~~~~~~~~~~~~~~~~~~~~~~~~~~~~~~~~~~~~~~~~~~(6.3)\newline
Now for $ 0\leq p \leq 1,$ using $(a + b)^{p} \leq a^{p} + b^{p},$ we can write
\begin{equation*}
{\mid t-x_{0}\mid }^{\alpha_{1}}\leq {\mid t-x\mid }^{\alpha_{1}} + {\mid x-x_{0}\mid }^{\alpha_{1}}
\end{equation*}
~~~~~~~~~~~~~~~~~and
\begin{equation*}
{\mid s-y_{0}\mid }^{\alpha_{2}}\leq {\mid s-y\mid }^{\alpha_{2}} + {\mid y-y_{0}\mid }^{\alpha_{2}}
\end{equation*}
Using these inequalities in(6.3), we get
\begin{equation*}
\mid K_{n_{1},n_{2}}(f;p_{n_{1}},q_{n_{1}};p_{n_{2}},q_{n_{2}}; x,y)-f(x,y)\mid~\leq~ K_{n_{1},n_{2}}({\mid t-x\mid }^{\alpha_{1}}{\mid s-y\mid }^{\alpha_{2}};p_{n_{1}},q_{n_{1}};p_{n_{2}},q_{n_{2}}; x,y)
\end{equation*}

~~~~~~~~~~~~~~~~~~~~~~~~~~~~~~~~~~~~~~~~~~~~~$+~{\mid y-y_{0} \mid}^{\alpha_{2}}K_{n_{1},n_{2}}({\mid t-x \mid}^{\alpha_{1}};p_{n_{1}},q_{n_{1}};p_{n_{2}},q_{n_{2}}; x,y)$\newline

~~~~~~~~~~~~~~~~~~~~~~~~~~~~~~~~~~~~~~~~~~~~~$+~{\mid x-x_{0} \mid}^{\alpha_{1}}K_{n_{1},n_{2}}({\mid s-y \mid}^{\alpha_{2}};p_{n_{1}},q_{n_{1}};p_{n_{2}},q_{n_{2}}; x,y)$\newline

~~~~~~~~~~~~~~~~~~~~~~~~~~~~~~~~~~~~~~~~~~~~~$+~{\mid x-x_{0} \mid}^{\alpha_{1}}{\mid y-y_{0} \mid}^{\alpha_{2}}K_{n_{1},n_{2}}(f_{0};p_{n_{1}},q_{n_{1}};p_{n_{2}},q_{n_{2}}; x,y)$

Now using the H¨older’s inequality for $ p_{1} = \frac{2}{\alpha_{1}}$,~~$p_{2} = \frac{2}{\alpha_{2}}$~~$q_{1} = \frac{2}{2-\alpha_{1}}$~~$q_{2} = \frac{2}{2-\alpha_{2}}$, we get

$K_{n_{1},n_{2}}({\mid t-x\mid }^{\alpha_{1}}{\mid s-y\mid }^{\alpha_{2}};p_{n_{1}},q_{n_{1}};p_{n_{2}},q_{n_{2}}; x,y)$\newline

$~=~K_{n_{1},n_{2}}({\mid t-x \mid}^{\alpha_{1}};p_{n_{1}},q_{n_{1}};p_{n_{2}},q_{n_{2}}; x,y)K_{n_{1},n_{2}}({\mid s-y \mid}^{\alpha_{2}};p_{n_{1}},q_{n_{1}};p_{n_{2}},q_{n_{2}}; x,y)$\newline

$~\leq~[(K_{n_{1},n_{2}}({\ t-x)^{2}};p_{n_{1}},q_{n_{1}};p_{n_{2}},q_{n_{2}}; x,y)]^{\frac{\alpha_{1}}{2}} {[K_{n_{1},n_{2}}(f_{0};p_{n_{1}},q_{n_{1}};p_{n_{2}},q_{n_{2}}; x,y)]}^{\frac{2-\alpha_{1}}{2}}$\newline

$~\times~[(K_{n_{1},n_{2}}({\ s-y)^{2}};p_{n_{1}},q_{n_{1}};p_{n_{2}},q_{n_{2}}; x,y)]^{\frac{\alpha_{2}}{2}} {[K_{n_{1},n_{2}}(f_{0};p_{n_{1}},q_{n_{1}};p_{n_{2}},q_{n_{2}}; x,y)]}^{\frac{2-\alpha_{2}}{2}}$\newline

This consequently gives the expected result. so the proof is complete.~~~~~~~~~~~~~$\square$\newline

\textbf{Remark 6.3.} If we take $ E = [0, \infty), $ then because of $ d(x, E) = 0$ and $ d(y, E) = 0$, we have
\begin{equation*}
\mid K_{n_{1},n_{2}}(f;p_{n_{1}},q_{n_{1}};p_{n_{2}},q_{n_{2}}; x,y)-f(x,y)\mid \leq M_({p_{n_{1}},q_{n_{1}},p_{n_{2}},q_{n_{2}}})^{4 -\frac{\alpha_{1}+ \alpha_{2}}{2}}{\delta_{n_{1}}(x)^{\frac{\alpha_{1}}{2}}}{\delta_{n_{2}}(y)^{\frac{\alpha_{2}}{2}}}
\end{equation*}

\textbf{Remark 6.4.} Using (5.3), it can be easily verified that $ st - lim_{n_{1}}{\delta_{n_{1}}} = 0 $ and $ st - lim_{n_{2}}{\delta_{n_{2}}} = 0. $
So we can estimate the order of statistical approximation of our bivariate operators by means of Lipschitz type maximal functions using this
result.

\end{document}